\documentclass{amsart}

\usepackage{amsfonts}
\usepackage{amsthm}
\usepackage{amsmath}
\usepackage{amsfonts}
\usepackage{latexsym}
\usepackage{amssymb}
\usepackage[latin1]{inputenc}
\usepackage{verbatim}

\newcommand{\U}{\mathcal{U}}
\newcommand{\N}{\mathcal{N}}

\newcommand{\um}{\U_{\cM}}

\def\um{\mathcal{U}_\m}
\newcommand{\uni}[1]{\um( #1 )}
\newcommand{\cluni}[1]{\overline{\um( #1 )}}
\newcommand{\cluno}[1]{\overline{\um( #1 )}^{\topologia}}

\def\a{\mathcal A}

\newcommand{\CC}{\mathbb{C}}
\newcommand{\RR}{\mathbb{R}}
\newcommand{\NN}{\mathbb{N}}

\newcommand{\D}{\mathcal D}

\newtheorem{fed}{Definition}[section]
\newtheorem{teo}[fed]{Theorem}
\newtheorem*{teo*}{Theorem}
\newtheorem{lem}[fed]{Lemma}
\newtheorem{cor}[fed]{Corollary}
\newtheorem{pro}[fed]{Proposition}
\theoremstyle{definition}

\newcommand{\topologia}{\,\mbox{\tiny $\sigma$-sot}}

\def\m{\mathcal{M}}
\def\msa{\m^{\text{sa}}}

\def\cA{\mathcal{A}}

\def\cD{\mathcal{D}}

\def\cN{\mathcal{N}}

\def\cZ{\mathcal{Z}}

\newcommand{\ds}{\displaystyle}

\newcommand{\hil}{\mathcal{H}}

\begin{document}

\title{A Schur-Horn theorem in II$_1$ factors}
\author{M. Argerami}
\address{Department of Mathematics, University of Regina, Regina SK, Canada}
\email{argerami@math.uregina.ca}
\thanks{M. Argerami supported in part by the Natural Sciences and Engineering Research
Council of Canada}
\author{P. Massey}
\address{Departamento de Matem\'atica, Universidad Nacional de La Plata and
Instituto Argentino de Matem\'atica-conicet, Argentina}
\email{massey@mate.unlp.edu.ar}
\thanks{P. Massey supported in part by CONICET of Argentina and a PIMS Postdoctoral Fellowship}



\begin{abstract} Given a  II$_1$ factor $\m$ and a diffuse abelian von Neumann subalgebra $\a\subset\m$,
we prove a version of the Schur-Horn theorem, namely
\[\overline{E_\a(\uni{b})}^{\topologia}=\{a\in \a^{sa}:\ a\prec b\},\ \ \ b\in\m^{sa},\]
where $\prec$ denotes spectral majorization, $E_\a$ the unique
trace-preserving conditional expectation onto $\a$, and $\U_\m(b)$
the unitary orbit of $b$ in $\m$. This result is inspired by a
recent problem posed by Arveson and Kadison. \end{abstract}

\maketitle

\section{Introduction}

In 1923, I. Schur \cite{schur} proved that if $A\in M_n(\CC)^{sa}$
(i.e., $A$ is selfadjoint) then
\[
\sum_{j=1}^k\,\alpha^\downarrow_j\le\sum_{j=1}^k\,\beta^\downarrow_j,\ k=1,\ldots,n,
\]
with equality when $k=n$ (denoted $\alpha\prec\beta$),
where $\alpha=\mbox{diag}(A)\in\RR^n$, $\beta=\lambda(A)\in\RR^n$ the spectrum
(counting multiplicity) of $A$,
and $\alpha^\downarrow,\beta^\downarrow\in\RR^n$ are obtained from $\alpha, \beta$
by reordering their entries in decreasing order.

In 1954, A. Horn \cite{Horn54} proved the converse: given $\alpha,\beta\in\RR^n$ with
$\alpha\prec\beta$, there exists a selfadjoint matrix $A\in M_n(\CC)$ such that $\mbox{diag}(A)=\alpha$,
$\lambda(A)=\beta$.
Since every selfadjoint matrix is diagonalizable, the
results of Schur and Horn can be combined in the following assertion: if $\cD$ denotes
the diagonal masa in $M_n(\CC)$ and $E_\cD$ is the compression onto $\cD$ , then
\begin{equation}\label{forma sh mat}
E_\cD(\{U\,M_\beta\,U^*: U\in M_n(\CC) \mbox{ unitary}\})=\{M_\alpha\in \cD:\ \alpha
\prec\beta\},
\end{equation}
where $M_\alpha$ is the diagonal matrix with the entries of
$\alpha$ in the main diagonal.

This combination of the two results, commonly known as Schur-Horn
theorem, has played a significant role in many contexts of matrix
analysis: although simple, vector majorization expresses a natural
and deep relation among vectors, and as such it has been a useful
tool both in pure and applied mathematics. We refer to the books
\cite{arnold,marshall} and the introductions of \cite{arvkad,neu}
for more on this.

During the last 25 years, several extensions of majorization have
been proposed by, among others, Ando \cite{Ando} (to selfadjoint
matrices), Kamei \cite{kam} (to selfadjoint operators in a II$_1$
factor), Hiai \cite{Hiai0,Hiai} (to normal operators in a von
Neumann algebra), and Neumann \cite{neu} (to vectors in
$\ell^\infty(\NN)$). With these generalizations at hand, it is
natural to ask about extensions of the Schur-Horn theorem.

In \cite{neu}, Neumann developed his extension of majorization
with the goal of using it to prove a Schur-Horn type theorem in
$\mathcal B(\hil)$ in the vein of previous works in convexity (see
the introduction in \cite{neu} for details and bibliography).
Other versions of the Schur-Horn theorem have
been considered in \cite{JMD} and \cite{arvkad}. It is interesting
to note that the motivation in \cite{neu} comes from geometry, in
\cite{JMD} comes from the study of frames on Hilbert spaces,
while in \cite{arvkad} it is of an operator theoretic nature.

In \cite{arvkad} Arveson and Kadison proposed the study of a
Schur-Horn type theorem in the context of II$_1$ factors, which
are for such purpose the most natural generalization of full
matrix algebras. They proved a Schur type theorem for II$_1$
factors and they posed as a problem a converse of this result,
i.e. a Horn type theorem. In this note we prove a Schur-Horn type
theorem that is inspired by Arveson-Kadison's conjecture (Theorem
\ref{la otra inclusion}).

\section{Preliminaries}\label{dos}

Throughout the paper $\m$ denotes a II$_1$ factor with normalized
faithful normal trace $\tau$. We denote by $\msa$, $\m^+$, $\um$,
the sets of selfadjoint, positive, and unitary elements of $\m$, and by $\cZ(\m)$ the center of $\m$. Given $a\in \m^{sa}$ we denote
its spectral measure by $p^a$.
The characteristic function of the set
$\Delta$ is denoted by $1_\Delta$. For $n\in \NN$,
the algebra of $n\times n$ matrices over
$\CC$ is denoted by $M_n(\mathbb C)$, and its unitary group by $\U_n$.
 By $dt$ we denote integration with respect to Lebesgue
measure. To simplify terminology, we will refer to non-decreasing
functions simply as ``increasing''; similarly, ``decreasing'' will
be used instead of ``non-increasing''.

Besides the usual operator norm in $\m$, we consider the
1-norm induced by the trace, $\|x\|_1=\tau(|x|)$. As we will
be always dealing with bounded sets in a II$_1$ factor, we can profit from the fact
that the topology induced by
${\|\cdot\|_1}$ agrees with the $\sigma$-strong operator topology.
Because of this we will express our results in terms of $\sigma$-strong closures
although our computations are based on estimates for $\|\cdot\|_1$.
For $X\subset \m$, we shall denote by $\overline X$ and
$\overline X^{\topologia}$ the respective closures in the norm
and in the $\sigma$-strong operator topology.

\subsection{Spectral scale and spectral preorders}

The
\emph{spectral scale} \cite{petz} of $a\in \m^{sa}$ is defined as
$$\lambda_a(t)=\min\{s\in \RR:\ \tau(p^{a}(s,\infty))\leq t\}, \ \
t\in [0,1).$$ The function $\lambda_a:[0,1)\rightarrow [0,\|a\|]$ is
decreasing and right-continuous. The map $a\mapsto\lambda_a$ is continuous with respect to both
$\|\cdot\|$ and $\|\cdot\|_1$, since \cite{petz}
\begin{equation}\label{dependencia cont 1}
\|\lambda_a-\lambda_b\|_\infty\leq \|a-b\|, \ \ \
\|\lambda_a-\lambda_b\|_1\leq \|a-b\|_1\ \ \ \
a,\,b\in\m^{sa},\end{equation} where the norms on the left are
those of $L^\infty([0,1],dt)$ and $L^1([0,1],dt)$ respectively.

\medskip

\noindent We say that $a$ is \emph{submajorized} by $b$, written
$a\prec_w b$, if $$ \int_0^s\lambda_a(t)\ dt\leq \int_0^s \lambda_b(t) \
dt,\ \ \ \text{for every }s\in [0,1). $$ If in addition
$\tau(a)=\tau(b)$ then we say that $a$ is \emph{majorized} by $b$,
written $a\prec b$.

\medskip
These preorders play an important role in many
papers (among them we mention \cite{fack, Hiai0, Hiai, kam}), and
they arise naturally in several contexts in operator theory and
operator algebras: some recent examples closely related to our
work are the study of Young's type \cite{Doug} and Jensen's type
inequalities \cite{JD,arvkad,silva,Hiai0}.

\begin{teo}[\cite{Hiai0}]\label{algo sobre domi} Let $a,\,b\in \m^{sa}$. Then
$a\prec b$ (resp. $a\prec_w b$) if and only if
$\tau(f(a))\leq \tau(f(b))$ for every
convex (resp. increasing convex) function $f:J\rightarrow \RR$, where $J$ is an open interval such that
$\sigma(a),\,\sigma(b)\subseteq J$.
\end{teo}

If $\N\subset \m$ is a von Neumann subalgebra and $b\in\m^{sa}$,
we denote by $\Omega_\N(b)$ the set of elements in $\N^{sa}$ that are majorized by $b$, i.e.
\[\Omega_\N(b)=\{a\in \N^{sa}:\ a\prec b\}.\ \]
The unitary orbit of $a \in\m^{sa}$ is the set $\U_\m(a)=\{ u^*au
: u \in \U_\m\}.$

\begin{pro}\label{la inclusion facil}
Let $\cN\subset \m$ be a von Neumann subalgebra  and let $E_\cN$
be the trace preserving conditional expectation onto $\cN$. Then,
for any $b\in\m^{sa}$,
\begin{enumerate}
\item\label{inclusion facil1} $E_\cN(b)\prec b$.
\item\label{inclusion facil2} $\|E_\cN(b)\|_1\leq \|b\|_1.$
\item\label{inclusion facil3}
$\overline{E_\cN(\U_\m(b))}^{\topologia}\subset \Omega_\N(b)$.
\end{enumerate}
\end{pro}
\begin{proof}(\ref{inclusion facil1}) The map $E_\cN$ is doubly stochastic (i.e. trace preserving, unital, and positive),
 so this follows from \cite[Theorem 4.7]{Hiai0}
(see also \cite[Theorem 7.2]{arvkad}).

(\ref{inclusion facil2}) Consider the convex function $f(x)=|x|$. Since $E_\cN(b)\prec b$, using
Theorem \ref{algo sobre domi} we get
$$\|E_\cN(b)\|_1=\tau(f(E_\cN(b)))\leq
\tau(f(b))=\|b\|_1.$$

(\ref{inclusion facil3})
 By (\ref{inclusion facil1}) and the fact that $ubu^*\prec b$ for every $u\in\U_\m$,
 we just have to prove that the set $\Omega_\N(b)$ is $\|\cdot\|_1$-closed. So, let $(a_n)_{n\in
\NN}\subset \Omega_\N(b)$ be such that $\lim_{n\rightarrow
\infty}\|a_n-a\|_1=0$ for some $a\in \cN$. Then, necessarily,
$a\in \cN^{sa}$. By (\ref{dependencia cont 1}), $$\int_0^s
\lambda_a(t)\ dt=\lim_{n\rightarrow
\infty}\int_0^s\lambda_{a_n}(t)\ dt\leq \int_0^s\lambda_b(t)\
dt.$$ Also, $\tau(a)=\lim_n\tau(a_n)=\tau(b)$, so $a\prec b$.
\end{proof}

The following result seems to be well-known, but we have not been able to find a reference.
Thus we give a sketch of a proof.
\begin{pro}\label{quien sabe que}
Let $a\in\cA^{sa}$, where $\cA\subset\m$ is a diffuse von Neumann subalgebra. Then there exists a spectral
resolution $\{e(t)\}_{t\in[0,1]}\subset\cA$ with $\tau(e(t))=t$ for every $t\in[0,1]$, and such that
\[
a=\int_0^1\,\lambda_a(t)\,de(t).
\]
\end{pro}
\begin{proof}
Since $\tau(1)<\infty$, it is enough to show the result for $a\ge0$. By \cite[Theorem 3.2]{Mas} there exists
$a'\in\cA$ with $g_{a'}(s)=\tau(p^{a'}(-\infty,s])$ continuous in $\RR$,
and an increasing left-continuous function $f$ such that $p^a(-\infty,s]=p^{a'}(-\infty,f(s)]$. Although the
original statement in \cite{Mas} involves a masa, only the fact that the algebra is diffuse is needed for its proof.

Let $g_{a'}^\dagger(t)=\min\{s:\,g_{a'}(s)\ge t\}$, and let $q(t)=p^{a'}(-\infty,g_{a'}^\dagger(t)]$. Since
$\tau(q(t))=g_{a'}(g_{a'}^\dagger(t))=t$ for every $t\in[0,1]$, it follows that $\{q(t)\}_{t\in[0,1]}$ is a
continuous spectral resolution. Moreover, $q(g_{a'}(t))=p^{a'}(-\infty,t]$, so
$p^a(-\infty,t]=q(g_{a'}(f(t)))$. As $g_{a'}\circ f$ is increasing and right-continuous, by \cite[Theorem 4.4]{Mas}
there
exists an increasing and left-continuous function $h_a$ such that $a=\int\,h_a(t)\,dq(t)$. Define $e(t)
=1-q(1-t)$, and $h(t)=h_a(1-t)$. Then $\tau(e(t))=t$, $a=\int h(t)\,de(t)$. As $h$ is decreasing and
right-continuous, it can be seen that $h=\lambda_a$.
\end{proof}

A spectral resolution $\{e(t)\}_{t\in[0,1]}$ as in Proposition \ref{quien sabe que} is called a {\em complete
flag} for $a$.

\section{A Schur-Horn Theorem for II$_1$ factors}\label{seccion sh}

For each $n\in\NN$, $k\in{1,\ldots,2^n}$, let $\{I^{(n)}_k\}_{k=1}^{2^n}$ be the partition of $[0,1]$
associated to the points $\{h\,2^{-n}:\ h=0,\ldots,2^n\}$.

\begin{fed}For each $n\in \NN$ and every $f\in
L^1([0,1])$, let
\begin{equation*}
E_n(f)= \sum_{i=1}^{2^n}\left( {2^n}\int_{I^{(n)}_i} f\
\right) \ 1_{I^{(n)}_i}. \end{equation*} \end{fed}

It is clear that each operator $E_n$ is a linear contraction for
both $\|\cdot\|_1$ and $\|\cdot\|_\infty$.

Given the flag $\{e(t)\}_t$, we write $e([t_0,t_1])$ for $e(t_1)-e(t_0)$. Note that, since we consider $e(t)$ diffuse,
$e([t_0,t_1])=e((t_0,t_1))=e([t_0,t_1))=e((t_0,t_1])$.

\begin{lem}\label{aproximaciones por discretizados}
Let $\{e(t)\}_{t\in[0,1]}$, $\{I^{(n)}_i\}_{i=1}^{2^n},\ n\in \NN$ and $\{E_n\}_{n\in \NN}$ as above. Then, for each $a\in\m^{sa}$,
\begin{equation}\label{ec lo que inter}
\lim_{n\rightarrow \infty}\|\,a-\int_0^1 E_n(\lambda_a)(t)\,de(t)\,\|_1=0.\end{equation}
\end{lem}

\begin{proof} By continuity of the trace, we only need to check that
\begin{equation}\label{ec lo que inter2}\|\lambda_a-E_n(\lambda_a)\|_1\xrightarrow[n]{}0\end{equation}
in $L^1([0,1])$. Consider first  a continuous function $g$. By uniform continuity, $\|g-E_n(g)\|_1\rightarrow0$.
Since continuous functions are dense in $L^1([0,1])$ and because the operators
$E_n$ are $\|\cdot\|_1$-contractive for every $n\in \NN$, a standard
$\varepsilon /3$ argument proves (\ref{ec lo que inter2}) for any integrable function.
\end{proof}

Recall that
$\cD$ denotes the diagonal masa in $M_n(\CC)$, and that for
$\alpha\in\RR^n$ we denote by $M_\alpha$ the matrix with the entries of
$\alpha$ in the diagonal and zero off-diagonal. The projection
$E_\cD$ of $M_n(\CC)$ onto $\cD$ is then given by
$E_\cD(A)=M_{\mbox{diag}(A)}$, where $\mbox{diag}(A)\in\RR^n$ is the main diagonal of $A$.
We use $\{e_{ij}\}$ to denote the
canonical system of matrix units in $M_n(\CC)$.

\medskip

\begin{lem}\label{dem} Let $\N\subset \m$ be a von Neumann
subalgebra,  and assume $E_\N$ denotes the unique trace preserving
conditional expectation onto $\N$. Let $\{p_i\}_{i=1}^n\subset
\cZ(\N)$ be a set of mutually orthogonal equivalent projections such that
$\sum_{i=1}^n p_i=I$.
Then there exists a unital *-monomorphism $\pi:M_n(\CC)\rightarrow
\m$ satisfying \begin{equation}\label{me pregunto si sera
necesario} \pi(e_{ii})=p_i, \ \ 1\leq i\leq n,
\end{equation}
\begin{equation}\label{comp con comp}
E_\N(\pi(A))=\pi(E_\D(A)),\ \ A\in M_n(\CC).\end{equation}
\end{lem}

\begin{proof} Since the projections $p_i$ are equivalent in $\m$, for each
$i$ there exists a partial isometry $v_{i1}^{}$ such that
$v_{i1}^{}v_{i1}^*=p_i$ and $v_{i1}^*v_{i1}^{}=p_1$. Let $v_{11}^{}=p_1$,
$v_{1i}^{}=v_{i1}^*$ for $2\leq i\leq n$ and $v_{ij}^{}=v_{i1}^{}v_{1j}^{}$
for $1\leq i,j\leq n$.
 In this way we get the standard associated system of matrix units $\{v_{ij}^{},\ 1\leq i,\,j\leq n\}$ in $\m$.
Define $\pi:M_n(\CC)\rightarrow \m$ by
$\pi(A)=\sum_{i,\,j=1}^n a_{ij}^{} \, v_{ij}^{}.$

The matrix unit relations imply that $\pi$ is a *-monomorphism and
it is clear that (\ref{me pregunto si sera necesario}) is also
satisfied. Moreover,
$$E_\N(v_{ij}^{})=E_\N(p_i\,v_{ij}^{}\,p_j)=p_{i}\,E_\N(v_{ij}^{})\,p_j=\delta_{ij}\,p_i,$$
since $p_i\,p_j=\delta_{ij}\,p_i$, $E_\N(v_{ii})=E_\N(p_i)=p_i$, and $p_i\in \cZ(\N)$.
Finally, we check (\ref{comp con comp}):
\[E_\N(\pi(A))=\sum_{i,j}a_{ij}E_\N(v_{ij})=\sum_ia_{ii}p_i=\pi(E_\D(A)).\qedhere\]
\end{proof}

Next we state and prove our version of the Schur-Horn
theorem for II$_1$ factors. Note the formal analogy with
(\ref{forma sh mat}).

\begin{teo}\label{la otra inclusion}
Let $\cA\subset \m$ be a diffuse abelian von Neumann subalgebra and let $b\in \m^{sa}$. Then
\begin{equation}\label{carac de la mayo esp}
\overline{E_\a(\um(b)) }^{\topologia}=\Omega_\cA(b).
\end{equation}
\end{teo}

\begin{proof}
By Proposition \ref{la inclusion facil}, we only need to prove
$\overline{E_\cA(\U_\m(b))}^{\topologia}\supset
\Omega_\cA(b)$.

So let $a\in \cA^{sa}$ with $a\prec b$. By Proposition \ref{quien sabe que},
\[a=\int_0^1\,\lambda_a(t)\,de(t),\ \ \ b=\int_0^1\,\lambda_b(t)\,df(t)\] where
$\{e(t)\}$ and $\{f(t)\}$ are complete flags with $\tau(e(t))=\tau(f(t))=t$, $e(t)\in\cA$, $t\in[0,1]$.

Let $\{I_i^{(n)}\}_{i=1}^{2^n}$, $n\in\NN$, be the family of partitions considered before
and let $\epsilon >0$. By Lemma \ref{aproximaciones por
discretizados} there exists $n\in \NN$ such that
\begin{equation}\label{los discretizados}
 \left\|a-\sum_{i=1}^{2^n} \alpha_i\ p_i\right\|_1< \epsilon, \ \ \
 \left\|b-\sum_{i=1}^{2^n} \beta_i\ q_i\right\|_1< \epsilon,
\end{equation} where $\alpha_i=2^n\,\int_{I_i^{(n)}}\,\lambda_a(t)\,dt$, $\beta_i=2^n\,\int_{I_i^{(n)}}\,\lambda_b(t)\,dt$,
$p_i=e(I_i^{(n)})$, $q_i=f(I_i^{(n)})$, $1\leq i\leq 2^n$. Note that $\tau(p_i)=\tau(q_i)=
2^{-n}$. Let $\alpha=(\alpha_1,\ldots,\alpha_{2^n})$, $\beta=(\beta_1,\ldots,\beta_{2^n})\in\RR^{2^n}$.
From the fact that $\lambda_a$ and $\lambda_b$ are decreasing, the entries of $\alpha$ and $\beta$ are already
in decreasing order. Using that $a\prec b$ IN $\m$, we conclude  that $\alpha\prec\beta$ in $\RR^n$.

By the classical Schur-Horn theorem \eqref{forma sh mat}, there
exists $U\in \U_n(\CC)$ such that
\begin{equation}\label{troesmas} E_\mathcal D(U M_{\beta}
U^*)=M_\alpha\end{equation} Consider the *-monomorphism $\pi$ of Lemma
\ref{dem} associated with the orthogonal family of projections
$\{p_i\}_{i=1}^{2^n}\subset \cA$. Let $w\in\um$ such that $wq_iw^*=p_i$, $i=1,\ldots,2^n$, and put
$u:=\pi(U)\,w\in \U_\m$. By
(\ref{comp con comp}) and (\ref{troesmas}),
\begin{eqnarray*}
\ds E_\cA\left(u\left(\,\sum_{i=1}^{2^n}\beta_i\,q_i\right)
u^*\right)&=& E_\cA(\pi(U M_\beta U^*))\\ \ds &=&
\pi(E_\mathcal D( U M_\beta U^*))=\pi(M_\alpha)\\ &=& \sum_{i=1}^{2^n} \alpha_i\,p_i
\end{eqnarray*} Using (\ref{los discretizados}),
we conclude that $$\left\|E_\a(ubu^*)-a\right\|_1\leq
\left\|E_\a(u(b-\sum_{i=1}^{2^n} \beta_i\
q_i)u^*)\right\|_1+\left\|a- \sum_{i=1}^{2^n} \alpha_i\
p_i\right\|_1< 2\epsilon.$$ As $\epsilon$ was arbitrary, we obtain
$a\in\overline{E_\cA(\U_\m(b))}^{\topologia}$.
\end{proof}

\begin{cor}\label{convexa y cerrada}
 For each $b\in\m^+$, the set $\overline{E_\a(\um(b))}^{\topologia}$ is convex and
$\sigma$-weakly compact.
\end{cor}

 In \cite{arvkad}, Arveson and
Kadison posed the problem whether for $b\in
\m^{sa}$, with the notations of Theorem \ref{la otra inclusion},
\begin{equation*}
E_\cA\left(\overline{\U_\m(b)}\right)=\Omega_\cA(b).
\end{equation*}
Since (\cite{kam}) $$\cluni{b}=\cluno{b}=\{a\in \m^{sa}:\ \lambda_a=\lambda_b\},$$
an affirmative
answer to the Arveson-Kadison problem is equivalent to
\begin{equation}\label{conj arv-kad}
E_\cA\left(\overline{\U_\m(b)}^{\topologia}\right)=\Omega_\cA(b).
\end{equation} As a description of the set $\Omega_\cA(b)$, (\ref{carac de la mayo esp})
is weaker than (\ref{conj arv-kad}), since in general
\begin{equation}\label{rela de esp}
E_\cA(\overline{\U_\m(b)}^{\topologia})\subset
\overline{E_\cA(\U_\m(b))}^{\topologia}.\end{equation}
An affirmative answer to the Arveson-Kadison problem would imply equality in (\ref{rela de esp}). We
think it is indeed the case, although a proof of this does not emerge from our present methods.

\noindent {\bf Acknowledgements.} We would like to thank
Professors D. Stojanoff, D. Farenick, and D. Sherman
for fruitful discussions regarding the material contained in this
note. We would also like to thank the referee for several suggestions that considerably simplified the exposition.

\end{document}